\documentclass[12pt]{article}
\usepackage{amssymb,latexsym}

\let\oldref\ref
\def\ref#1{{\normalfont\oldref{#1}}}
\def\eqref#1{{\normalfont(\oldref{#1})}}

\def\qedsymbol{\vbox{\hrule\hbox{%
\vrule height1.3ex\hskip0.8ex\vrule}\hrule}}

\def\myendproof{\qquad\qedsymbol}
\def\noqed{\def\qedsymbol{}}

\def\N{\mathbb{N}}

\def\Nset{\mathbf{N}}
\newtheorem{thm}{Theorem}

\newtheorem{lem}{Lemma}

\title{A Bilateral Version of the Shannon-McMillan-Breiman Theorem}

\author{%
Pierre Tisseur\\
Laboratoire G\'enome et Informatique\\
Universit\'e d'Evry, Tour Evry 2.}

\begin{document}

\maketitle

\begin{abstract}
We give a new version of the Shannon-McMillan-Breiman theorem in the
case of a bijective action.
For a finite partition $\alpha$ of a compact set $X$ and a
measurable action $T$ on $X$, we denote by
$C_{n,m,\alpha }^{T}(x)$ the element of the partition
$\alpha\vee T^{1}\alpha\vee \ldots \vee T^{m}\alpha\vee T^{-1}\alpha
\vee \ldots \vee T^{-n}\alpha$ which contains a point $x$.
We prove that for $\mu$-almost all $x$,
$$
\lim_{n+m\to\infty}\left(\frac{-1}{n+m}\right)
\log\mu (C_{n,m,\alpha}^T(x))
=h_\mu (T,\alpha ),
$$
where $\mu$ is a $T$-ergodic probability measure and
$h_\mu (T,\alpha)$ is the
metric entropy of $T$ with respect to the partition $\alpha$.

\end{abstract}


\section{Introduction}

The Shannon-McMillan-Breiman theorem \cite{brei57},
\cite{pete83} is used in many problems
related to the metric entropy map of an ergodic measure.
We extend this well-known result to the
case of a bijective dynamical system.
Our proof follows the line of
Petersen's proof \cite{pete83}.
We illustrate this new result with an example
that gives an inequality between shifts and cellular automata
entropies and some analog of the Lyapunov exponents.
Our bilateral version of the Shannon-McMillan-Breiman theorem is expected
to be useful in other areas of dynamical systems.

\section{Background material}

Let $X$ be a compact space,
$\mu$ a probability measure on $X$
and $T$ a measurable map from $X$ to $X$.
We denote by $\alpha$ a finite partition of $X$ and by
$C_{n,\alpha }^{T}(x)$ the element of the partition
$\alpha\vee T^{-1}\alpha\vee \ldots \vee
T^{-n}\alpha$ which contains the point $x$.
For all point $x$ the information map $I$ is defined by
$$
I(\alpha )(x) = -\log \mu(C_{\alpha}^T(x))
=\sum_{A\in\alpha}-\log \mu (A)\chi_A(x),
$$
where $C_{\alpha}^T(x)$ is the element of $\alpha$ which contains $x$
and
$\chi_A$ is the characteristic function defined by
$$
\chi_A(x) = \left\{\begin{array}{l}
\mbox {$1$ if $x\in A$},\\
\mbox {$0$ otherwise.}
\end{array}\right.
$$
The information map satisfies
\begin{eqnarray}
I(\alpha\vee \beta )&=&
I(\alpha )+I(\beta \vert\alpha ), \label{P1}\\
I(T\alpha \vert T\beta )&=&I(\alpha \vert \beta )\circ T^{-1}
\label{P2}.
\end{eqnarray}
These two properties are easily proved from the definition of $I$
and the fact that $T$ is a surjective map.
We refer to
\cite[p. 238]{pete83}, \cite[Chap.8]{pomi98}
for a detailed proof of \eqref{P1} and \eqref{P2}.

A simple formulation of the metric entropy with respect to the partition $\alpha$ is given by
$$
h_\mu(T,\alpha)=\lim_{n\to\infty}\int_X I(\alpha \vert
\vee_{k=1}^n T^k \alpha)(x)\,d\mu (x),
$$
where
$$
I(\alpha \vert\beta )(x)=
-\sum_{A\in\alpha, B\in\beta} \chi_{A\cap B}(x)
\log\left(\frac{\mu(A\cap B)}{\mu (B)}\right)
$$
is the conditional information map
representing the quantity of information given by the partition
$\alpha$ knowing the partition $\beta$ about the point $x$.


We recall the Shannon-McMillan-Breiman theorem \cite{brei57}
\cite{pete83}.

\begin{thm}[Shannon-McMillan-Breiman's theorem]
If $\mu$ is a $T$-ergodic measure, then
for $\mu$-almost all $x$ in a compact $X$ we have
$$
\lim_{n\to\infty}\frac{-1}{n}\log\mu (C_{n,\alpha}^T(x))=h_\mu (T,\alpha).
$$
\end{thm}
\section{A bilateral version of Shannon-McMillan-Breiman's theorem}

In order to prove the main result (Theorem \ref{th4}) we need to expose two technical lemmas.
The proof of Lemma \ref{lem3} and Theorem \ref{th4} requires
a bilateral version of the Birkhoff pointwise ergodic theorem:
for a $T$-ergodic measure $\mu$ one has
$$
\lim_{n+m\to\infty}\frac{1}{n+m+1}\sum_{k=-m}^n f\circ T^k (x)
=\int_X f(x)d\mu (x)
$$
for almost all $x$ with a map $f$ in $L_1$.
This result is easily deduced from the Birkhoff pointwise ergodic
theorem (see \cite[Chap.10]{pomi98}) by
breaking up the infinite sum in two proportional parts.

\begin{lem}\label{lem2}
For all integers $m$ and $n$ we have
$$
\sum_{k=-m}^{n-1} I(\alpha \vert \vee_{j=1}^{n-k}T^{-j}\alpha )\circ
T^{k}=I(\vee_{j=-n}^m T^j \alpha )-I(T^{-n}\alpha ).
$$
\end{lem}
 
\begin{proof}
Note that $I(\vee_{j=-n}^m T^j \alpha)= I(\vee_{j=0}^{m+n}T^{m-j}\alpha)$.
Using \eqref{P1} we have
$$
I(\vee_{j=0}^{m+n}T^{m-j}\alpha)=I(\vee_{j=1}^{m+n}T^{m-j}\alpha)
+I(T^m\alpha \vert \vee_{j=1}^{m+n}T^{m-j}\alpha)
$$
and from \eqref{P2} we get
$$
I(T^m\alpha \vert \vee_{j=1}^{m+n}T^{m-j}\alpha)=
I(\alpha \vert \vee_{j=1}^{m+n}T^{-j}\alpha )\circ T^{-m}.
$$
Hence,
$$
I(\vee_{j=-n}^m T^j \alpha)=
I(\vee_{j=1}^{m+n}T^{m-j}\alpha )+
I(\alpha \vert \vee_{j=1}^{m+n}T^{-j}\alpha )\circ T^{-m}.
$$
The same operations on
$I(\vee_{j=1}^{m+n}T^{m-j}\alpha)$
yields
\begin{eqnarray*}
I(\vee_{j=1}^{m+n}T^{m-j}\alpha )
&=&I(\vee_{j=2}^{m+n}T^{m-j}\alpha )+I(\alpha \vert \vee_{j=2}^{m+n}T^{1-j}\alpha ) \circ T^{-m+1}\\
&=&I(\vee_{j=1}^{m+n-1}T^{m-1-j}\alpha )+I(\alpha
\vee_{j=1}^{m+n-1}T^{-j}\alpha ) \circ T^{-m+1}.
\end{eqnarray*}
Hence,
\begin{eqnarray*}
I(\vee_{j=-n}^m T^j \alpha)
&=&I(\vee_{j=1}^{m+n-1}T^{m-1-j}\alpha )\\
&&{}+I(\alpha\vee_{j=1}^{m+n-1}T^{-j}\alpha ) \circ T^{-m+1}+
I(\alpha \vert \vee_{j=1}^{m+n}T^{-j}\alpha )\circ T^{-m}.
\end{eqnarray*}
 
Iterating similarly $t-1$ times on $I(\vee_{j=1}^{m+n-1}T^{m-1-j}\alpha )$
leads to
$$
I(\vee_{j=-n}^m T^j \alpha )= I(\vee_{j=1}^{m+n-t}T^{m-t-j}\alpha )+\sum_{k=0}^{t} I(\alpha \vert \vee_{j=1}^{m+n-k}T^{-j}\alpha )\circ T^{-m+k}.
$$
Taking $t=m+n-1$ gives
$$
I(\vee_{j=-n}^m T^j \alpha )= +I(T^{-n}\alpha )+\sum_{k=0}^{m+n-1} I(\alpha \vert \vee_{j=1}^{m+n-k}T^{-j}\alpha )\circ T^{-m+k}
$$
which completes the proof.

\end{proof}

\begin{lem}\label{lem3}
If $\mu$ is a $T$ ergodic measure then
for almost all $x$ in $X$,
\begin{eqnarray*}
\lim_{m+n\to\infty}\frac{1}{m+n+1}\sum_{k=-m}^{n-1}\lim_{s\to\infty}
I(\alpha \vert \vee_{j=1}^{s}T^{-j}\alpha )\circ T^k(x)&=&\\
&& \hspace{-6cm}
\lim_{m+n\to\infty}\frac{1}{m+n+1}\sum_{k=-m}^{n-1}I(\alpha \vert
\vee_{j=1}^{n-k} T^{-j}\alpha )\circ T^k(x).
\end{eqnarray*}
\end{lem}

\begin{proof}
For notational convenience, we introduce
$$
f=\lim_{s\to\infty}I(\alpha \vert\vee_{j=1}^{s}T^{-j}\alpha)
\quad \mbox{and} \quad
F_N=\sup_{s\ge N}\vert
I(\alpha \vert\vee_{j=1}^{s}T^{-j}\alpha )-f\vert.
$$
It is well known that the sequence
$(I(\alpha\vert\vee_{j=1}^{s}T^{-j}\alpha))_{s\in\N}$ converge
almost everywhere and in $L_1$.
The proof of this convergence (see \cite[Chap.8]{pomi98} and \cite[p.262]{pete83})
 requires the increasing martingale theorem.

We need to show that
$$
\lim_{m+n\to\infty}\frac{1}{m+n+1}
\sum_{k=-m}^{n-1}\left \vert I(\alpha \vert
\vee_{j=1}^{n-k}T^{-j}\alpha )\circ T^k-f \circ T^k\right\vert =0.
$$
Note that
\begin{eqnarray*}
\frac{1}{m+n+1}\sum_{k=-m}^{n-1}\vert
I(\alpha \vert \vee_{j=1}^{n-k}T^{-j}\alpha )
\circ T^k-f\circ T^k\vert
&\le&\\
&& \hspace{-5cm}
\frac{1}{m+n+1}\sum_{k=n-N}^{n-1}\vert
I(\alpha \vert \vee_{j=1}^{n-k}T^{-j}\alpha )
\circ T^k-f\circ T^k\vert \\
&&\hspace{-5cm}+
\frac{1}{m+n+1}\sum_{k=-m}^{n-N-1}\vert F_N\vert\circ T^k .
\end{eqnarray*}
If we fix $N$ and let $n+m$ tend to infinity then the first term in
the right-hand side of the above inequality goes to zero.
Since the map $F_N$ belongs to $L_1$ (see
\cite{pete83}), the bilateral version of Birkhoff's ergodic theorem
applies and we can assert that
$$
\lim_{m+n\to\infty}\frac{1}{m+n+1}\sum_{k=-m}^{n-N-1}\vert F_N\vert\circ T^k = \int_X F_N \, d\mu .
$$
Since $\lim_{N\to\infty}F_N=0$,
the dominated convergence theorem implies that
$\int_X F_N \, d\mu$ tends to zero which completes the proof.

\end{proof}

\begin{thm}\label{th4}
For a bijective map $T$ from $X$ to $X$ and a $T$-ergodic measure
$\mu$, we have for $\mu$-almost all $x$
$$
\lim_{n+m\to\infty}\frac{-1}{n+m}\log\mu (C_{n,m,\alpha}^T(x))
=h_\mu (T,\alpha ),
$$
where $C_{n,m,\alpha}^T(x)$ represents the element of the partition
$\alpha\vee T\alpha \ldots \vee T^m\alpha\vee T^{-1}\alpha\vee \ldots \vee
T^{-n}\alpha $ containing the point $x$.
\end{thm}

\begin{proof}
Since the sequence
$\left (I(\alpha \vert\vee_{j=1}^{s}T^{-j}\alpha )\right )_{s\in\N}$
converges to a $L_1$ map
and by using the dominated convergence theorem, it follows that
\small
$$
h_\mu (T,\alpha )
=\lim_{s\to\infty}\int_XI(\alpha \vert \vee_{j=1}^{s}T^{-j}
\alpha )(x)\,d_\mu (x)
=\int_X\lim_{s\to\infty}I(\alpha \vert\vee_{j=1}^{s}T^{-j}
\alpha )(x)\,d_\mu (x)
$$
\normalsize
for $\mu$-almost all $x$.
The bilateral version of Birkhoff's ergodic theorem implies that for almost all $x$
$$
h_\mu (T,\alpha )=\lim_{m+n\to\infty}\frac{1}{m+n+1}\sum_{k=-m}^{n-1}\lim_{s\to\infty}
I(\alpha \vert \vee_{j=1}^{s}T^{-j}\alpha )\circ T^k(x).
$$
{}From Lemma \ref{lem3} it follows that
$$ h_\mu (T,\alpha ) = \lim_{m+n\to\infty}\frac{1}{m+n+1}\sum_{k=-m}^{n-1}I(\alpha \vert\vee_{j=1}^{n-k}
T^{-j}\alpha )\circ T^k(x).
$$
Using Lemma \ref{lem2} for almost all $x$ we obtain
$$
h_\mu (T,\alpha ) = \lim_{m+n\to\infty}\frac{1}{m+n+1}\left( I(\vee_{j=-n}^mT^j\alpha )(x)-I(T^{-n}\alpha)(x)
\right).
$$
Since
$(I(T^{-n}\alpha))_{n\in\Nset}$ is bounded for $\mu$-almost all point $x$, the sequence
$\frac{I(T^{-n}\alpha)}{n+m+1}$ tends almost surely to zero.
Hence,
\begin{eqnarray*}
h_\mu (T,\alpha )&=&
\lim_{m+n\to\infty}\frac{1}{m+n+1}
I(\vee_{j=-n}^m T^j \alpha)(x)\\
&=&
-\lim_{m+n\to\infty}\frac{1}{m+n}
\log\mu (C_{n,m,\alpha}^T(x)).
\myendproof
\end{eqnarray*}

\noqed
\end{proof}

\section{An illustration}

In this example we do not give a definition of the particular discrets dynamical
systems called cellular automata;
the reader can find a survey in \cite{bkm97}
and the complete proof of this illustration in \cite{tiss99}. The bilateral
version of the Shannon-McMillan-Breiman theorem is needed to establish an inequality
between the entropy of a cellular automaton $F$ denoted by $h_\mu(F,\alpha )$,
the entropy $h_\mu (\sigma )$ of a particular bijective cellular automaton
$\sigma$ called the shift,
and some discret analog of the Lyapunov exponents.
Here the
measure $\mu$ is $F$-invariant and $\sigma$-ergodic.

The standart Shannon-McMillan-Breiman theorem
\cite[Chap.10]{pomi98} says that in the case
of an invariant measure $\mu$.
\[
h_\mu (F,\alpha ) = \int \lim_{n\to\infty}\frac{-1}{n}\log \mu \left (C^F_{n,\alpha} (x)\right ).
\]
In \cite{tiss99} one proves that there exists some integer and bounded maps
$f_n$ and $g_n$ such that, for all point $x$, one has
$$
C^F_{n,\alpha }(x)\supset C^\sigma_{f_n(x),g_n(x),\alpha}(x)
$$
with
$\lim_{n\to\infty}f_n(x)+g_n(x)=+\infty$ for $\mu$-almost all point $x$ for a
certain class of cellular automata. For those that do not belong to this class,
the entropy is equal to zero (see \cite{tiss99}).
With these properties we obtain
\[
h_\mu (F,\alpha ) \le \int \lim_{n\to\infty}\frac{-1}{n}
\log \mu \left (C^\sigma_{f_n(x),g_n(x),\alpha} (x)\right )d\mu (x)
\]
and 
\begin{eqnarray*}
h_\mu (F,\alpha ) \le \int\liminf_{n\to\infty}\frac{-1}{f_n(x)+g_n(x)+1}
\log \mu \left (C^\sigma_{f_n(x),g_n(x),\alpha} (x)\right )\\
\times
 \frac{g_n(x)+f_n(x)+1}{n}d\mu (x).
\end{eqnarray*}
The bilateral version of the Shannon-McMillan-Breiman theorem implies that 

\[
h_\mu (F,\alpha )\le h_\mu (\sigma,\alpha )\int \liminf_{n\to\infty}\frac{f_n(x)+g_n(x)+1}{n}.
\]
Using the Fatou lemma, we have
\[
h_\mu (F,\alpha ) \le h_\mu (\sigma, \alpha )\times (\lambda^+_\mu +\lambda^-_\mu),
\]
where
$$
\lambda^+_\mu = \liminf\int \frac{f_n(x)}{n}d\mu (x)
\quad \mbox{and} \quad
\lambda^-_\mu =
\liminf\int \frac{g_n(x)}{n}d\mu (x)
$$
are called the left and right average Lyapunov exponents.

\bibliographystyle{siam}

\end{document}